\documentclass[10pt,twocolumn,twoside]{IEEEtran}
\IEEEoverridecommandlockouts

\usepackage[cmex10]{amsmath}
\usepackage{color}
\usepackage{amssymb}
\usepackage{algorithmic}
\usepackage{color}
\usepackage[titlenumbered,ruled]{algorithm2e}
\usepackage{dsfont}
\usepackage{balance}
\usepackage{graphicx}
\usepackage{tikz}
\usepackage{subcaption}
\usepackage{balance}
\usepackage{graphicx}

\newtheorem{definition}{Definition}
\newtheorem{theorem}{Theorem}

\newtheorem{lemma}{Lemma}

\newcommand{\scalemath}[2]{\scalebox{#1}{\mbox{\ensuremath{\displaystyle #2}}}}

\newcommand{\llqs}{\ensuremath{\lambda_l(\mathbf{Q}_S)}}
\newcommand{\ls}{\ensuremath{\lambda_l}}

\newcommand{\tr}[1]{\ensuremath{\textbf{tr}\left(#1\right)}}

\newcommand{\N}{\mathcal{N}}

\newcommand{\Da}{\mathbf{D}_\kappa}
\newcommand{\daa}{d_{\kappa,a}}
\newcommand{\vaa}{\kappa_a}

\newcommand{\Q}[1]{\ensuremath{\mathbf{Q}_{#1}}}

\newcommand{\Qit}[1]{\ensuremath{\mathbf{Q}_{#1}^{-2}}}
\newcommand{\Qi}[1]{\ensuremath{\mathbf{Q}_{#1}^{-1}}}

\newcommand{\A}{\ensuremath{\mathbf{A}}}
\newcommand{\B}{\ensuremath{\mathbf{B}}}
\newcommand{\C}{\ensuremath{\mathbf{C}}}
\newcommand{\x}{\ensuremath{\mathbf{x}}}
\newcommand{\umat}{\ensuremath{\mathbf{u}}}
\newcommand{\w}{\ensuremath{\boldsymbol{\mu}}}
\newcommand{\LL}{\ensuremath{\mathbf{L}}}
\newcommand{\Ds}{\ensuremath{\mathbf{D}_S}}
\newcommand{\D}[1]{\ensuremath{\mathbf{D}_{#1}}}
\newcommand{\I}{\ensuremath{{\mathbf{I}}}}
\newcommand{\W}{\ensuremath{\mathbf{W}}}
\newcommand{\Z}{\ensuremath{\mathbf{Z}}}
\newcommand{\F}{\ensuremath{\mathbf{F}}}
\newcommand{\G}{\ensuremath{\mathbf{G}}}
\newcommand{\J}{\ensuremath{\mathbf{J}}}
\newcommand{\K}{\ensuremath{\mathbf{K}}}
\newcommand{\On}{\ensuremath{{\mathbf{0}}_n}}
\newcommand{\Om}{\ensuremath{\mathbf{0}_m}}

\newcommand{\bQ}[1]{\ensuremath{b_1\mathbf{Q}_{#1}}}

\newcommand{\bQi}[1]{\ensuremath{\big(b_1\mathbf{Q}_{#1}\big)^{-1}}}
\newcommand{\bbQi}[1]{\ensuremath{\big(b_2\mathbf{Q}_{#1}-b_3\bf{I}\big)^{-1}}}

\title{\LARGE \bf
Submodularity in Systems with Higher Order Consensus with Absolute Information
}

\author{Erika Mackin and Stacy Patterson
    \thanks{This work was funded in part by NSF awards CNS-1527287 and CNS-1553340.} 
\thanks{ Erika Mackin and Stacy Patterson are with the Department of Computer Science, Rensselaer Polytechnic Institute, Troy, New York 12180, USA. Email:
    {\tt\small mackie2@rpi.edu, sep@cs.rpi.edu} } }

\begin{document}

\maketitle
\thispagestyle{empty}
\pagestyle{empty}

\begin{abstract}
We investigate the performance of $m^{th}$ order consensus systems with stochastic external perturbations, where a subset of leader nodes incorporates absolute information into their control laws.  
The system performance is measured by its coherence, an $\mathcal{H}_2$ norm that quantifies the total steady-state variance of the deviation from the desired trajectory.  
We first give conditions under which such systems are stable, and we derive expressions for coherence in stable second, third, and fourth order systems. We next study  the problem of how to identify a set of leaders that optimizes coherence. To address this problem, we define  set functions that quantify each system's coherence and prove that these functions are submodular. This allows the use of an efficient greedy  algorithm that to find a leader set with which coherence is within a constant bound of optimal. We demonstrate the performance of the greedy algorithm empirically, and further, we show that the optimal leader sets for the different orders of consensus dynamics do not necessarily coincide.
\end{abstract}


%

\section{Introduction}
Consensus algorithms and their applications have long been an important field of research. 
One of the most common real-world applications of consensus dynamics is the control of autonomous vehicles~\cite{BJMP12, HB13, HBV10}. This application has been studied widely in terms of first and second order dynamics, in which vehicles compare the differences between their position and velocity and that of their neighbors to keep themselves in a desired formation. Consensus dynamics that incorporate vehicles' measurements of their acceleration and jerk have been considered as well, but not yet as widely. By including these higher order dynamics in the control mechanisms, the group of vehicles is able to respond to abrupt changes in direction and speed, which requires adjustments to their acceleration, and potentially other higher order states as well~\cite{RMC06}. These sudden changes may be caused by rough terrain, high winds, unexpected obstacles, etc. \cite{RMC06, BSSV15}. Other applications of higher order consensus dynamics include unmanned aerial~\cite{DYSZ15} or underwater vehicles~\cite{WYL12} and swarms of satellites~\cite{SSL07}.

It has been shown that, by incorporating absolute state into the control law of every vehicle, the formation can more closely follow the desired trajectory than can a formation that uses dynamics based solely on the relative measurements of the difference between each vehicles' state and the states of its neighbors~\cite{BJMP12}. In real-world applications, the vehicles in a formation may not all have the computational power and energy resources to support hardware that generates these absolute measurements. For example, they may be located in a remote area that makes physical access for upgrades difficult, or it may simply be a matter of a limited budget that does not allow for purchasing a new GPS for each vehicle. Further, even in systems with homogeneous nodes, there may be cost or communication limitations that prohibit the use of such hardware in every vehicle. Given these limitations, it is important to judiciously select which agents have the capability to access absolute information to optimize the performance of the formation. 

When only some vehicles have access to absolute information, the group is using \emph{leader-follower consensus}, a subset of consensus dynamics 
where  \emph{leader} nodes update their states using both absolute and relative information, and the remaining nodes, the \emph{followers} use only relative information in their control laws. 
Leader-follower dynamics has been widely studied, albeit primarily in first and second order dynamics~\cite{PB10, LFJ11, BJMP12, HBV10}. We consider higher order consensus dynamics with stochastic external perturbations entering through the highest-order state variable only. The performance of the system is measured by its \emph{coherence}, the total steady-state deviation of each node's first order state from its desired value.
%
 
 We pose the problem of selecting the set of leaders that minimizes the coherence of the formation in second, third, and fourth order systems. The coherence can be computed only when the system is stable. We first derive the conditions under which stability is guaranteed for each order. We also prove that systems of order four and higher are always unstable when all gains are identical. We then derive the closed-form expression for coherence for stable second, third, and fourth order systems.
We prove that all three coherence expressions can be reformulated as submodular set functions; the submodularity property can be informally thought of as a property of diminishing returns as set size increases.
Since all three set functions are submodular, we can use a computationally-inexpensive greedy algorithm to identify a set of leader nodes, such that this set yields coherence that is within a provable bound of the coherence of the optimal leader set~\cite{NWF78}.  We further demonstrate that the performance of the greedy algorithm often yields a leader set that has coherence very close to that of the optimal set. 
We also show that the optimal leader sets for the different orders of consensus dynamics do not necessarily coincide.
 

\subsubsection*{Related work} 
The problem of leader selection for optimal coherence has been studied widely in first order systems \cite{PSFS17, PB10, LFJ11, FL16, PMD16, P17, LFJ14, DCJ16, PYZ17}. 
Of significant note to this paper are the recent works \cite{CBP14} and \cite{MP17-2}, which show that the leader selection problem can be expressed as a submodular optimization problem.
The first work considers noise-free leaders, while the second considers noise-corrupted leaders, as we do in this work.

Second order leader-follower consensus systems have been studied, though not yet as widely. In \cite{BJMP12}, the coherence of a system where all nodes have access to absolute information  is studied for various graph topologies. The convergence requirements for a system with a single noise-free leader is studied in \cite{QLMC14}. Systems with multiple noise-free leaders are studied in terms of controllability conditions in \cite{GR10}. The stability margin in systems where the communication graph is a lattice and leaders make up one or more edges of the lattice is studied in \cite{HBV10}. 

Systems with third order and higher dynamics have received less attention to date. The problem of higher order consensus with a leader was first introduced in \cite{RMC06}, where the authors proved that consensus is reachable in third order systems only when the system is stable. Consensus in third order systems was studied in vehicle platoons with a single leader \cite{BSSV15} and in directed networks without leaders in \cite{HZX18}.
Consensus without leaders in higher order systems was studied in \cite{JWJ09}, where the presence of a spanning tree in the network topology was proved to be a necessary condition for consensus, 
and in \cite{RA15}, where the problem of developing a control input capable of guiding all nodes to the same location was studied. 
Finally, in \cite{TBS19} the authors prove that consensus is not guaranteed to be reachable in higher order systems as the network grows infinitely large. Further, they  
give conditions for consensus in third order systems where all nodes are leaders.
To the best of our knowledge, no previous work has studied coherence in systems of order three and higher. 
 
 In a preliminary version of this work~\cite{MP18-2}, we studied the leader selection problem in second order systems.
 We studied only the case where both gains on the system were equal to one. We derived the expression for coherence and proved that a set function based on the coherence was submodular.
 In this paper, we expand on these results by allowing the gains in the second order system to be any values that satisfy the stability requirements. We also provide stability conditions for third and fourth order systems, and we prove that set functions based on the coherence of third and fourth order systems are  submodular. 

In the remainder of the paper, we first describe our system model and problem formulation in Section~\ref{model.sec}. We present stability conditions and derive the expression for coherence in second, third, and fourth order systems in Section \ref{stabcoh.sec}. We give some background on submodularity and then prove that the set functions on coherence for second, third, and fourth order systems are submodular in Section~\ref{submod.sec}. We present our numerical results in Section~\ref{numAn.sec}, before concluding in Section~\ref{concl.sec}.

\section{System Model} \label{model.sec}
We study consensus dynamics in a connected, undirected network, which is modeled by a graph $\mathcal{G}=(V,E,W)$, where $V$ is the set of $n$ nodes, $E$ is the set of edges, and ${W: E \rightarrow \mathds{R}}$ is a weight function that assigns a positive value to each edge in $E$. 

Let each node $i$ have $m$ scalar-valued states, denoted by $x_{j}^{i}$, $j=1,\ldots m$, where $\x_j$ is the $n$-vector of order $j$ node states.  A subset of nodes are selected to be leaders, which then receive external information on each of their $m$ states. All other nodes update their states using only the state values transmitted by their neighbors. 

Without loss of generality, we assume that the desired formation is for all nodes to converge to position $0$.
When the system is of order $m$, then the nodes update their states as follows:
\begin{align*}
\dot{\x}_j &=\x_{j+1}, \: \: j=1, \ldots, m-1\\
 \dot{\x}_m &= \umat + \w
\end{align*} 
where $\w$ is an $n$-vector of zero-mean, white stochastic disturbances and $\umat$ is an $n$-vector of the nodes' control inputs. The control input for node $i$ is defined as:
\[ u_i = - \sum_{j=1}^{m} a_j \left(\sum_{k \in \N_i} w_{ik}(x_{j}^{i}-x_{j}^{k}) +\delta_i\kappa_i x_{j}^{i} \right) \]
where $x_{j}^{i}$ is the $i^{th}$ element of vector $\x_j$, $a_j$ is a non-zero fixed gain, $\N_i$ is the set of neighbors of node $i$, $w_{ik}$ is the weight of edge $(i, k)$, $\kappa_i$ is a positive scalar, and $\delta_i$ is an indicator variable with value $1$ when node $i$ is a leader and $0$ otherwise. 

Let $\Da$ be a diagonal matrix where $\Da(i,i)=\kappa_i >0$. Let $\Ds$ be a diagonal matrix corresponding to the set of leaders $S$, where $\Ds(i,i)=\delta_i$. We define the matrix $\mathbf{Q}_S$ as ${\mathbf{Q}_S=\LL+\Da \Ds}$, where $\LL$ is the weighted Laplacian matrix of $\mathcal{G}$, such that
 \begin{align*}
 \LL_{jk} =
 \begin{cases}
      -w_{jk} & (j,k) \in E  \\
      \sum_{i=1}^n w_{ji} & j=k \\
      0 & \text{otherwise}.
\end{cases}
\end{align*} 

When the set $S$ consists of a single node, so that $S=\{j\}$, we denote the corresponding matrix as $\mathbf{Q}_j$.
We define the state vector of the entire system $\x \in \mathbb{R}^{n\times m}$ as
\[
\x= \begin{bmatrix}
    \x_1^T & \x_2^T   &\ldots & \x_{m-1}^T  & \x_m^T
\end{bmatrix} ^T .
\]
The states of all nodes are then updated as: 
\begin{align*}
\dot{\x}
=
\A
\x
+
\B\w,
\end{align*}
where 
\begin{align}
\A&=
\begin{bmatrix}
   \On & \I_n & \On & \hdots & \On\\
   \On & \On & \I_n & \hdots & \On\\
   \vdots &\vdots & \vdots & \ddots & \vdots\\
   \On & \On & \On & \hdots & \I_n \\
  -a_1\Q{S} & -a_2\Q{S} & -a_3\Q{S} & \hdots & -a_m\Q{S}
\end{bmatrix} \label{A.eq},\\
\B &=  \begin{bmatrix}
    \On  & \On & \ldots & \On & \I_n 
\end{bmatrix}^T, \nonumber 
\end{align}
where $\On$ is an $n \times n$ matrix of all zeros and $\I_n$ is the $n \times n$ identity matrix. 

The output of the system, $\bf{y}\in\mathds{R}^n$, is studied in terms of first order states only, so that
\begin{align}
\bf{y} &= \C \x, \nonumber 
\end{align}
with
\begin{align*}
\C &=  \begin{bmatrix}
    \I_n  & \On & \ldots & \On  
\end{bmatrix}.
\end{align*}
The performance of the system, for a given set of leaders $S$ and weight matrix $\Da$, is quantified by the total steady-state variance of the deviations of the nodes' first order states from the desired trajectory. This variance is defined as: 
\[
H(S) =   \lim_{t \rightarrow \infty} \sum_{i=1}^{n} \textbf{var}(y_i) = \lim_{t \rightarrow \infty}  \sum_{j=1}^{n} \textbf{var}(x_{1}^{j}).
\]
This measure, also known as coherence, is bounded when the system is stable.
Further, for a stable system, the variance can be derived from the expression of the square of the system's ${\cal H}_2$ norm
\begin{align}
H(S) =  \tr{ \int_0^{\infty}\C e^{t\A} \B \B^T e^{t\A^T}\C^T dt},\nonumber 
\end{align}
or, alternatively,
\begin{align}
H(S)=\tr{\C \textbf{P}\C^T}, \nonumber 
\end{align}
where $\textbf{P}$ is the controllability Gramian, which is the solution to the Lyapunov equation
\begin{align}
\A \textbf{P}+\textbf{P}\A^T+\B\B^T=0. \label{lyap.eq}
\end{align}

When the system is stable and the expression for coherence is known, we can then turn to the problem of selecting the leader set that minimizes the coherence. We define the \emph{$m^{th}$ order $k$-leader selection problem}, for a non-negative integer $k$, as follows:
\begin{equation}
\begin{array}{ll}
\text{minimize} & H_m(S) \\
\text{subject to} & |S| \leq k.
\end{array} \label{Hprob.eq}
\end{equation}

\section{Stability and Coherence in Higher Order Systems} \label{stabcoh.sec}
We first prove under what conditions systems of order $m$ are stable and then derive expressions for the coherence of the stable systems.

\subsection{Stability Analysis}
The system is stable only when $Re(\lambda_i(\A))<0$, ${i=1, \ldots, n}$. 
To study the eigenvalues of $\A$, we decompose the matrix into
$n$ $m \times m$ matrices $\A_l$, ${l=1, \ldots, n}$, that can be analyzed in terms of the eigenvalues of $\Q{S}$.

Since ${\bf Q}_S$ is symmetric positive definite for any $S \neq \emptyset$ and any values $\kappa_i > 0$~\cite{ME10}, it can be diagonalized by a unitary matrix $\bf{U}$ whose columns are the eigenvectors of $\Q{S}$. 
Let $\bf{\Lambda}$ be the corresponding diagonal matrix of eigenvalues, so that ${\bf{\Lambda}=\bf{U}^{T}\Q{S} \bf{U}}$.
Define $\bf{\mathcal{U}}$ as the $nm \times nm$  block diagonal matrix, with each diagonal block equal to ${\bf{U}}$. Then,
\begin{align}
\bf{\mathcal{U}^{T}}\A\bf{\mathcal{U}}= 
\begin{bmatrix}
\On & \I_n & \On & \hdots & \On \\
\On & \On & \I_n & \hdots & \On\\
\vdots & \vdots &  \vdots & \ddots & \vdots \\
\On & \On & \On & \hdots & \I_n \\
-a_1\bf{\Lambda} & -a_2\bf{\Lambda} & -a_3\bf{\Lambda} & \hdots & -a_m\bf{\Lambda}
\end{bmatrix}. \label{blockeigmatrix}
\end{align}

Note that we can permute the rows and columns of (\ref{blockeigmatrix}) to obtain a block diagonal matrix of the form
\begin{align*}
\mathbf{P}\bf{\mathcal{U}}^{T}\A\bf{\mathcal{U}}\mathbf{P}=
\begin{bmatrix}
\A_1 & \Om & \hdots & \Om \\
\Om & \A_2  & \hdots & \Om\\
\vdots &  \vdots & \ddots & \vdots \\
\Om & \Om &  \hdots & \A_n
\end{bmatrix}
\end{align*}
where $\mathbf{P}$ is a permutation matrix. Each $m \times m$ matrix ${\bf A}_{l}$, ${l=1, \ldots, n}$, is of the form
\begin{align}
\A_{l}&= 
\begin{bmatrix}
0 & 1 &  \hdots & 0 \\
\vdots &   \vdots & \ddots & \vdots \\
0 & 0 & \hdots & 1 \\
-a_1\lambda_l(\Q{S}) & -a_2\lambda_l(\Q{S}) &  \hdots & -a_m\lambda_l(\Q{S})
\end{bmatrix}. \nonumber 
\end{align}

We note that the eigenvalues of $\A$  are the union of the eigenvalues of each $\A_l$, and thus $\A$ is stable if and only if $Re(\lambda_i(\A_{l}))<0$ for $l=1, \ldots n$, $i=1, \ldots, m$. 
To find when $Re(\lambda_i(\A_{l}))$ is guaranteed to be negative, we consider the $m \times m$ Hurwitz matrix $\textbf{M}$ of the characteristic polynomial of $\A_l$. The Hurwitz determinants, the determinants of the $m$ leading principal submatrices of $\textbf{M}$, are denoted by $\Delta_j$, ${j=1,\ldots,m}$. When all $\Delta_j$ are all positive, then $Re(\lambda_i(\A_{l}))<0$ for all $i=1, \ldots, m$, ${l=1, \ldots, n}$, and $\A$ is thus stable \cite{H95}. For simplicity, in the remainder of this section, we use $\lambda_l$ to denote $\lambda_l(\Q{S})$. Recall that ${\lambda_l=\lambda_l(\LL + \Da \Ds)}$ and the effect of the values $\kappa_i$ is included in the eigenvalues of $\A$. We, therefore, restrict our study of the necessary conditions for $\lambda_l<0$ to the values $a_i$ only.

\subsubsection{Second Order Stability}
It has been previously proven that $\A$ is stable when $m=2$ and $\Ds=\I$ \cite{PB10}. For completeness, we prove here that $\A$ is also stable when $\Ds \not = \I$.
\begin{theorem}
When $m=2$, $\A$ is stable if and only if the gains $a_1, \: a_2$ are positive. 
\end{theorem} \label{stable2.thm}
\begin{IEEEproof}
$\A$ is stable when all eigenvalues of $\A_{l}$ have negative real parts, for $l=1, \ldots n$, which is true if and only if the Hurwitz determinants of the characteristic polynomial of $\A_l$,  $p_{l}(s)=s^2+a_2\ls s + a_1\ls$, are positive. The Hurwitz matrix is 
\begin{align*}
\textbf{M}&=\begin{bmatrix} 
a_2\ls & 0 \\
1&  a_1\ls 
\end{bmatrix}
\end{align*}
and so we find that the two Hurwitz determinants of $\textbf{M}$ are 
\begin{align*}
\Delta_1&=a_2\ls\\
\Delta_2&=a_1a_2\ls^2.
\end{align*}

Thus, we can see that $\A$ is stable if and only if $a_1, a_2> 0$.
\end{IEEEproof}

\subsubsection{Third Order Stability} 
We now prove when $\A$ is stable for $m=3$.

\begin{theorem}  \label{stable3.thm}
When $m=3$, $\A$ is stable if and only if the gains are such that $a_1, a_2, a_3 >0$ and $\frac{a_2a_3}{a_1}\ls > 1$.
\end{theorem}
\begin{IEEEproof}
$\A$ is stable 
if and only if the Hurwitz determinants of the characteristic polynomial of $\A_l$, $p_{l}(s)=s^3+a_3\ls s^2+a_2\ls s + a_1\ls$, are positive. 
The Hurwitz matrix is 
\begin{align*}
\textbf{M}&=\begin{bmatrix} 
 a_3\ls & a_1 \ls & 0\\
1&  a_2\ls & 0 \\
0&  a_3\ls & a_1\ls
\end{bmatrix}
\end{align*}
and so we find that the three Hurwitz determinants of $\textbf{M}$ are 
\begin{align}
\Delta_1&=a_3\ls  \label{D1.eq} \\
\Delta_2&=a_2a_3\ls^2-a_1\ls \label{D2.eq} \\
\Delta_3&= a_1a_2a_3\ls^3-a_1^2\ls^2. \label{D3.eq}
\end{align}
%
We can see that (\ref{D1.eq}) is positive if and only $a_3\ls > 0$. It is straightforward to note that (\ref{D2.eq}) and (\ref{D3.eq}) cannot both be positive unless all gains are positive and $\frac{a_2a_3}{a_1} \ls> 1$,
which concludes the proof.
\end{IEEEproof}

\subsubsection{Higher Order Stability}
We now prove when $\A$ is stable for $m=4$ and $\A$ defined as in (\ref{A.eq}). We also prove that when $m\geq4$ and all $a_i$ are set to be equal, that $\A$ is never stable. 
\begin{theorem} \label{stable4.thm}
When $m=4$, $\A$ is stable if and only if the gains are such that $a_1, a_2, a_3, a_4>0$, $\frac{a_3a_4}{a_2}\ls > 1$, and ${\left(\frac{a_3a_4}{a_2}-\frac{a_1a_4^2}{a_2^2}\right) \ls>1}$.
\end{theorem}
\begin{IEEEproof}
$\A$ is stable 
if and only if the Hurwitz determinants of the characteristic polynomial of $\A_l$,  $p_{l}(s)=s^4+a_4\ls s^3+a_3\ls s^2+a_2\ls s + a_1\ls$, are positive. The Hurwitz matrix is 
\begin{align*}
\textbf{M}&=\begin{bmatrix} 
 a_4\ls & a_2 \ls & 0 & 0\\
1&  a_3\ls & a_1\ls & 0 \\
0&  a_4\ls & a_2\ls & 0 \\
0 & 1 & a_3 \ls & a_1 \ls
\end{bmatrix}
\end{align*}
and so we find that the four Hurwitz determinants of $\textbf{M}$ are 
\begin{align*}
\Delta_1&=a_4\ls  \\
\Delta_2&=a_3a_4\ls^2-a_2\ls\\
\Delta_3&= a_2a_3a_4\ls^3-a_1a_4^2 \ls^3 -a_2^2\ls^2 \\
\Delta_4&=a_1a_2a_3a_4\ls^4-a_1^2a_4^2 \ls^4 -a_1a_2^2\ls^3.
\end{align*}
We first note that $\Delta_1>0$ only when $a_4>0$. It is straightforward to observe, further, that if any of $a_1, a_2, a_3, a_4$ are negative, then some $\Delta_i$ is non-positive. We can see that $\A$ has all negative eigenvalues only when ${\frac{a_3a_4}{a_2}\ls>1}$, and ${(\frac{a_3a_4}{a_2}-\frac{a_1a_4^2}{a_2^2}) \ls>1}$, for all $l=1, \ldots, n$,  thus concluding the proof.
\end{IEEEproof}

In second and third order systems, it is possible for the stability conditions to be satisfied when all gains are equal. In fourth order systems, however, the system is never stable when all gains are the same. In fact, we find that systems of order $m \geq 4$ are never stable when all gains are equal.
\begin{theorem}
When $m\geq 4$ and $a_i=a$ for $i=1, \ldots, m$, $\A$ is not stable. 
\end{theorem}
\begin{IEEEproof}
First, consider the first and third Hurwitz determinants of the characteristic polynomial ${p_{l}(s)=s^m+a\ls s^{m-1}+\ldots+a\ls s + a\ls}$ for $\A_{m}$, when $m=4$ and $a>0$. They are:
\begin{align*}
\Delta_1&= a\ls \\
\Delta_3 &=
-a^2\ls^2. 
\end{align*}
We know that the eigenvalues of $\Q{S}$ are all real and so it can never be the case that both $a\ls>0$ and $-(a\ls)^2>0$. Therefore $\A$ is not stable when $m=4$ and all $a_i=a$. 

Next,  consider the third Hurwitz determinant when $m >4$ and $a>0$, 
\begin{align*}
\Delta_3 &=
 \begin{vmatrix} 
 a\ls & a\ls & a\ls\\
1&  a\ls & a\ls \\
0&  a\ls & a\ls
\end{vmatrix}=-a^2\ls^2 +a^2\ls=0.
\end{align*}
Since $\Delta_3$ is $0$ and thus non-positive for all $m>4$, $\A$ is not stable. 

When $a\leq0$, $\Delta_1$ is never positive for any ${m\geq 4}$. Therefore, for $m\geq4$ and $a_i=a$, $\A$ is never stable. 
\end{IEEEproof}

\subsection{Coherence Analysis} \label{coh.sec}
Now that we have proved under which conditions second, third, and fourth order systems are stable, we can derive expressions for their coherence.

\subsubsection{Second Order Coherence}

\begin{theorem} 
When $m=2$ and $\A$ is stable, the coherence of the system is:
\begin{align*}H_2(S)= \frac{1}{2a_1a_2}\tr{\Qit{S}}.\end{align*} 
\end{theorem}
\begin{IEEEproof}
By inspection, we suppose $\textbf{P}$ to be:
\begin{align*}
\textbf{P}=
\begin{bmatrix}
\frac{1}{2 a_1a_2}\Qit{S} & \On \\
\On & \frac{1}{2 a_2}\Qi{S}
\end{bmatrix}.
\end{align*}
It is straightforward to verify that $\mathbf{P}$ is the solution to the Lyapunov equation (\ref{lyap.eq}).
Therefore,  $ H_2(S)= \tr{{\bf C} {\bf P} \bf{C}^T} = \frac{1}{2 a_1a_2}\tr{\Qit{S}}$.
\end{IEEEproof}

\subsubsection{Third Order Coherence}
\begin{theorem}
When $m=3$ and $\A$ is stable, the coherence of the system is:
\begin{align*}H_3(S)= \frac{a_3}{2 a_1^2}\tr{\Qi{S}\left(\frac{a_2a_3}{a_1}\Q{S}-\I \right)^{-1}}. \end{align*}
\end{theorem}
\begin{IEEEproof}
By inspection, we suppose $\textbf{P}$ to be:
\begin{align}
\textbf{P}=
\begin{bmatrix}
\F & \On & -\G\\
\On & \G & \On \\
-\G & \On & \J
\end{bmatrix}, \label{P3}
\end{align}
where 
\begin{align*}
\F&=\frac{a_3}{2a_1}\Qi{S}(a_2a_3\Q{S}-a_1\I)^{-1} \\
\G&=\frac{1}{2}\Qi{S}(a_2a_3\Q{S}-a_1\I)^{-1}\\
\J&=\frac{a_2}{2}(a_2a_3\Q{S}-a_1\I)^{-1}.
\end{align*}
We first note that $a_1\F=a_3\G$, $a_2\Q{S}\G=\J$, and that $a_1\Q{S}\G-a_3\Q{S}\J=-\frac{1}{2}\I$, and therefore, when we substitute (\ref{P3}) into the Lyapunov equation (\ref{lyap.eq}), the expression holds.
Thus, \[H_3(S)= \frac{a_3}{2 a_1^2}\tr{\Qi{S}\left(\frac{a_2a_3}{a_1}\Q{S}-\I \right)^{-1}}.\]
\end{IEEEproof}

Note that when all gains have the same value $a$ the expressions for the coherence of second and third order systems can be written as 
\begin{align*}
H_2(S)&=\sum_{l=1}^n \frac{1}{(a\lambda_i(\Q{S}))^2} \\
H_3(S)&=\sum_{l=1}^n \frac{1}{a\lambda_i(\Q{S}) \left( a\lambda_i(\Q{S})-1\right)}.
\end{align*} Since both systems are stable, by Theorem \ref{stable3.thm}, $a \lambda_i(\Q{S}) >1$ for $i=1, \ldots, n$. From this, we can clearly see that \begin{align}{H_3(S) > H_2(S)} \label{O23coh.eq} \end{align} for all $S$.
\subsubsection{Fourth Order Coherence}
\begin{theorem}
When $m=4$ and $\A$ is stable, the coherence of the system is:
\begin{align*}
&H_4(S)= \frac{1}{2 a_1a_2}\textbf{tr}\Bigg(\Qit{S}\left(\frac{a_3a_4}{a_2}\Q{S} -\I\right)\\
&~~~~\left(\left(\frac{a_3a_4}{a_2}-\frac{a_1a_4^2}{a_2^2}\right)\Q{S}-\I\right)^{-1}\Bigg).
\end{align*}
\end{theorem}
\begin{IEEEproof}
By inspection, we suppose $\textbf{P}$ to be:
\begin{align}
&\textbf{P}=
\begin{bmatrix}
\F & \On & -\G & \On \\
\On & \G & \On & -\J \\
-\G & \On & \J & \On\\
\On & -\J & \On & \K
\end{bmatrix}, \label{P4}
\end{align}
where 
\begin{align*}
\F&=\scalemath{0.95}{\frac{1}{2 a_1}\Qit{S}\left(\frac{a_3a_4}{a_2}\Q{S} -I\right)\left(\left(a_3a_4-\frac{a_1a_4^2}{a_2}\right)\Q{S}-a_2\I\right)^{-1}}\\
\G&=\frac{a_4}{2 a_2}\Qi{S}\left(\left(a_3a_4-\frac{a_1a_4^2}{a_2}\right)\Q{S}-a_2\I\right)^{-1}\\
\J &=\frac{1}{2}\Qi{S}\left(\left(a_3a_4-\frac{a_1a_4^2}{a_2}\right)\Q{S}-a_2\I\right)^{-1}\\
\K &=\frac{1}{2}\left(a_3-\frac{a_1a_4}{a_2}\right)\left(\left(a_3a_4-\frac{a_1a_4^2}{a_2}\right)\Q{S}-a_2\I\right)^{-1}.
\end{align*}
We note that the following hold:
\begin{align*}
a_2\G&=a_4\J \\
 \J&=-a_1\Q{S}\F+a_3\Q{S}\G \\
 \K&=-a_1\Q{S}\G+a_3\Q{S}\J \\
 -\frac{1}{2}\I&=a_2\Q{S}\J-a_4\Q{S}\K.
\end{align*}

  Thus, when we substitute (\ref{P4}) into the Lyapunov equation (\ref{lyap.eq}) the expression holds and, therefore,
\begin{align*}
&H_4(S)= \frac{1}{2 a_1a_2}\textbf{tr}\Bigg(\Qit{S}(\frac{a_3a_4}{a_2}\Q{S} -\I)\\
&~~~~\left(\left(\frac{a_3a_4}{a_2}-\frac{a_1a_4^2}{a_2^2}\right)\Q{S}-\I\right)^{-1}\Bigg).
\end{align*}
\end{IEEEproof}

We next use the coherence expressions to study the $k$-leader selection problem for second, third, and fourth order systems.


\section{Leveraging Submodularity for Leader Selection} \label{submod.sec}

As an alternative to the leader selection problem (\ref{Hprob.eq}) defined in Section \ref{model.sec}, we can also define an optimization problem of the form
\begin{equation} \label{general.prob}
\begin{array}{ll}
\text{maximize} & f_m(S) \\
\text{subject to} & |S| \leq k, 
\end{array} 
\end{equation}
where $f_m$ is the set function $f_m: 2^V \rightarrow  \mathbf{R}$
\begin{equation*} 
 f_m(S) = \left\{ \begin{array}{ll} 
0 & \text{if}~S = \emptyset \\
 C_m - \rho_mH_m(S)& \text{otherwise}. 
\end{array} \right.
\end{equation*}
Let $\rho_m$ be a positive scalar and ${C_m = 2 \left(\max_{ s \in V}  \rho_mH_m(s)\right)}$. 
We note that maximizing $f_m(S)$ is equivalent to minimizing $H_m(S)$.

While a problem like (\ref{general.prob}) can be solved by an exhaustive search of all subsets of $V$ of size less than or equal to $k$, this approach becomes computationally infeasible for anything but small values of $k$ as the size of the network increases. 
Instead, we can approximate the optimal set of leaders with a \emph{greedy algorithm}, a more computationally tractable
 approach. The greedy algorithm is as follows: the  set of leaders $S$ is initialized to the empty set. In each iteration, the node ${v \in V \setminus S}$, which, when added to $S$, maximizes the value of $f_m$, is identified and added to $S$. After $k$ rounds, or when no further improvement is possible, the algorithm terminates. 

If $f_m$ is a non-decreasing, submodular set function, then the greedy algorithm generates a solution that is within a constant factor of optimal.
We first provide some background on submodularity and its applicability to greedy algorithms and then present our proofs for the submodularity of $f_m$ when $m=2,3,4$.

\vspace{-.3cm}
\subsection{Background} 
We make use of the following definitions and theorems to prove the submodularity of $f_m(S)$.
\begin{definition}[\cite{NWF78}]
A  function $f:2^V \mapsto\mathbf{R}$, where $V$ is a finite set, is called \emph{submodular} if, for all $A, B \subseteq V$,
\[ 
f(A) + f(B) \geq f(A \cup B) - f(A \cap B).
\]
\end{definition}

\begin{definition}
A set function $f:2^V \mapsto \mathbf{R}$ is called \emph{non-increasing} if for all  $A,B \subseteq V$,
if $A \subseteq B$, then $f(A) \geq f(B)$. The function $f$ is called \emph{non-decreasing} if for all $A, B \subseteq V$,
if $A \subseteq B$, then $f(A) \leq f(B)$.
\end{definition}

\begin{theorem}[\cite{NWF78} Prop. 4.3]  \label{submodBound.thm}
Let  $f:2^V \rightarrow \mathbf{R}$ be a non-decreasing, submodular set function. Let $f^*$ be the value of $f$ for an optimal set $S$ of size $k$ and let $S^g$ be the set of size $k$ returned by the greedy algorithm. 
Then,
\begin{align}
\frac{f^\star -f(S^g)}{f^* -f(\emptyset)} \leq \left( \frac{k-1}{k}\right) ^k \leq \frac{1}{e}. \label{submodBound.eq}
\end{align}
\end{theorem}


\subsection{Submodularity Analysis}




We now show that the set functions $f_m$, $m=2,3,4$, are submodular.
\subsubsection{Second Order Submodularity}

For systems where ${m=2}$, we define the set function 
\begin{equation*} 
 f_2(S) = \left\{ \begin{array}{ll} 
0 & \text{if}~S = \emptyset \\
 C_2 -\rho_2H_2(S) & \text{otherwise}, 
\end{array} \right.
\end{equation*}
where ${C_2=2 \left(\max_{ s \in V}  \rho_2H_2(s)\right)}$ and ${\rho_2 = 2a_1a_2}$. 
The function $f_2(S)$ is maximized when $H_2(S)$ is minimized.

\begin{theorem} \label{submod2.thm}
The set function $f_2(S)$ is non-decreasing and submodular.
\end{theorem}
The proof of this theorem depends on the following lemma, whose proof is reserved to the appendix.
\begin{lemma} \label{O23.lem}
Consider a set function $f$ defined as
\begin{equation} 
 f(S) = \left\{ \begin{array}{ll} 
0 & \text{if}~S = \emptyset \\
 C - \tr{\bQi{S}\bbQi{S}} & \text{otherwise}, \label{setf.eq}
\end{array} \right.
\end{equation}
where ${C = 2 \left(\max_{ s \in V}  \tr{\bQi{s}\bbQi{S}}\right)}$, $b_1, b_2>0$, $b_3 \geq 0$, $b_2\lambda_l(\Q{S})>b_3$ for $l=1, \ldots, n$, and both $\bQ{S}$ and $b_2\Q{S}-b_3\I$ are nonsingular. 
The function $f$ is non-decreasing and submodular. 
\end{lemma}

The proof for Theorem \ref{submod2.thm} is as follows.

\begin{IEEEproof}
We first note that ${\rho_2H_2(S) = \tr{\Qit{S}}}$ and $\Q{S}$ is non-singular, by definition.
Let $b_1=b_2=1$ and let $b_3=0$ so that (\ref{setf.eq}) is equal to $f_2$. We can then see that $f_2$ is non-decreasing and submodular. 
\end{IEEEproof}

\subsubsection{Third Order Submodularity}
For systems where $m=3$, we define the set function
\begin{equation*} 
 f_3(S) = \left\{ \begin{array}{ll} 
0 & \text{if}~S = \emptyset \\
 C_3 - \rho_3H_3(S) & \text{otherwise}, 
\end{array} \right.
\end{equation*}
where ${C_3=2 \left(\max_{ s \in V}  \rho_3H_3(s)\right)}$ and ${\rho_3=\frac{2a_1^2}{a_3}}$. The function $f_3(S)$ is maximized when $H_3(S)$ is minimized.

\begin{theorem} \label{submod3.thm}
The set function $f_3(S)$ is non-decreasing and submodular.
\end{theorem}
\begin{IEEEproof}
We first note that \[\rho_3H_3(S)=\tr{\Qi{S}\left(\frac{a_2a_3}{a_1}\Q{S}-\I\right)^{-1}}.\] We also note that, by Theorem \ref{stable3.thm}, $\frac{a_2a_3}{a_1}\llqs>1$ for ${l=1, \ldots, n}$, and, therefore, $\Q{S}$ and $\frac{a_2a_3}{a_1}\Q{S}-\I$ are both non-singular matrices.
Let $b_1=1$, $b_2=\frac{a_2a_3}{a_1}$, and $b_3=1$ so that (\ref{setf.eq}) is equal to $f_3$. We can then see that $f_3$ is non-decreasing and submodular. 
\end{IEEEproof}

\subsubsection{Fourth Order Submodularity}
For systems where $m=4$, we define the set function
\begin{equation*} 
 f_4(S) = \left\{ \begin{array}{ll} 
&0  ~~~\text{if}~S = \emptyset \\
 &C_4 - \rho_4H_4(S)~\text{otherwise}, 
\end{array} \right.
\end{equation*}
where 
$C_4 = 2 \left(\max_{ s \in V} \rho_4H_4(s) \right)$
and $\rho_4=2a_1a_2$.
When $f_4(S)$ is maximized, $H_4(S)$ is minimized.

\begin{theorem} \label{submod4.thm}
The set function $f_4(S)$ is non-decreasing and submodular.
\end{theorem}
The proof of this theorem depends on the following lemma, whose proof similar to that of Lemma~\ref{O23.lem} and is deferred to a technical report \cite{MP19}.

\begin{lemma} \label{O4.lem}
Consider a set function $f$ defined as
\begin{equation} 
 f(S) = \left\{ \begin{array}{ll} 
&0 ~~~~ \text{if}~S = \emptyset \\
 &C - \textbf{tr}\Big(\Qit{S}(b_1\Q{S}-\I)\\
 &~~~\big((b_1-b_2)\Q{S}-\I\big)^{-1}\Big)  ~\text{otherwise}, \label{setf4.eq}
\end{array} \right.
\end{equation}
where \[C = 2 \left(\max_{ s \in V}  \tr{\Qit{s}(b_1\Q{s}-\I)((b_1-b_2)\Q{s}-\I)^{-1}}\right),\] $b_1> b_2>0$, $(b_1-b_2)\lambda_l(\Q{S})>1$ for $l=1, \ldots, n$, and $\Q{S}$ and $(b_1-b_2)\Q{S}-\I$ are nonsingular. The function $f$ is non-decreasing and submodular. 
\end{lemma}


\begin{IEEEproof}
We first note that 
\begin{align*}\rho_4H_4(S)&=\textbf{tr}\Bigg(\Qit{S}\left(\frac{a_3a_4}{a_2}\Q{S} -\I\right)\\
&~~~~\left(\left(\frac{a_3a_4}{a_2}-\frac{a_1a_4^2}{a_2^2}\right)\Q{S}-\I\right)^{-1}\Bigg). \end{align*}
We also note that, by Theorem \ref{stable4.thm}, ${(\frac{a_3a_4}{a_2}-\frac{a1_a4^2}{a_2^2})\llqs>1}$ for $l=1, \ldots, n$, and, therefore, $\Q{S}$ and ${(\frac{a_3a_4}{a_2}-\frac{a1_a4^2}{a_2^2})\Q{S}-I}$ are both non-singular matrices.
Let ${b_1=\frac{a_3a_4}{a_2}}$ and ${b_2=\frac{a_1a_4^2}{a_2^2}}$ so that (\ref{setf4.eq}) is equal to $f_4$. We can then see that $f_4$ is non-decreasing and submodular. 
\end{IEEEproof}

\subsection{Performance of the Greedy Algorithm}

Now that we have proved that $f_2$, $f_3$, and $f_4$  are all non-decreasing and submodular, we can apply Theorem \ref{submodBound.thm} to bound the performance of the greedy algorithm for $k$-leader selection for each of these objective functions.
The following theorem follows directly from Theorems \ref{submodBound.thm}, \ref{submod2.thm}, \ref{submod3.thm}, and \ref{submod4.thm}.

\begin{theorem} \label{cohBound.thm}
Let $m \in \{2, 3, 4\}$.
Let $S_m^g$ be the set of $k$ leaders output by the greedy algorithm, and let $S^*_m$ be a set of $k$ leaders that minimizes $H_m(S)$. Then,
\begin{align*}
H_m(S_m^g) \leq \frac{C_m}{\rho_me} +\Big(1-\frac{1}{e}\Big)H_m(S^*_m).
\end{align*}
\end{theorem}

Finally, we analyze the computational complexity of the greedy algorithm.

When $m=2$, the complexity of the greedy algorithm is $O(kn^3)$. We first compute $\LL^\dagger$, an $O(n^3)$ operation. In each iteration $j$ of the greedy algorithm, we compute the rank-one update, an $O(n^2)$ operation, then find the trace of the square of the resulting matrix, which is also $O(n^2)$. We do this $|V|-j-1$ times; therefore, each iteration $j$ has complexity $O((|V|-j-1)n^2)$,
and there are at most $k$ iterations.

When $m=3$, the complexity of the greedy algorithm is again $O(kn^3)$. We first compute $\LL^\dagger$ and $(\frac{a_2a_3}{a_1}\LL-\I)^{-1}$, both of which are $O(n^3)$ operations. In each iteration $j$,  for each node $v \in V \setminus S_{j-1}$ we compute a rank-one update to find both $\Qi{S_{j-1} \cup \{v\}}$ and $\left(\frac{a_2a_3}{a_1}\Q{S_{j-1} \cup \{v\}} -\I \right)^{-1}$, an $O(n^2)$ operation for both matrices. We then find the trace of the product of the two matrices, which also has complexity $O(n^2)$. Thus, each iteration has complexity $O(3(|V|-j-1)n^2)$, and again, there are at most $k$ iterations.

When $m=4$, the expression for the coherence can be rearranged to be:
\begin{align*}
&H_4(S)=\\
&~~\textbf{tr}\Bigg(\Qit{S} - \frac{a_1a_4^2}{a_2^2}\Qi{S}\left(\left(\frac{a_3a_4}{a_2}-\frac{a_1a_4^2}{a_2^2}\right)\Q{S}-\I\right)^{-1}\Bigg).
\end{align*}
The complexity of the greedy algorithm  can then be analyzed in a similar manner to the previous two cases, and is again $O(kn^3)$.

\section{Experimental Results} \label{numAn.sec}

In this section, we first compare the optimal coherence of first, second, and third order systems. The set function based on the coherence of a first order system is submodular and the coherence of a first order system is defined as $H_1(S)=\frac{1}{2}\frac{1}{a_1}\tr{\Qi{S}}$~\cite{MP17-2}. 


\begin{figure}

 \includegraphics[scale=.44]{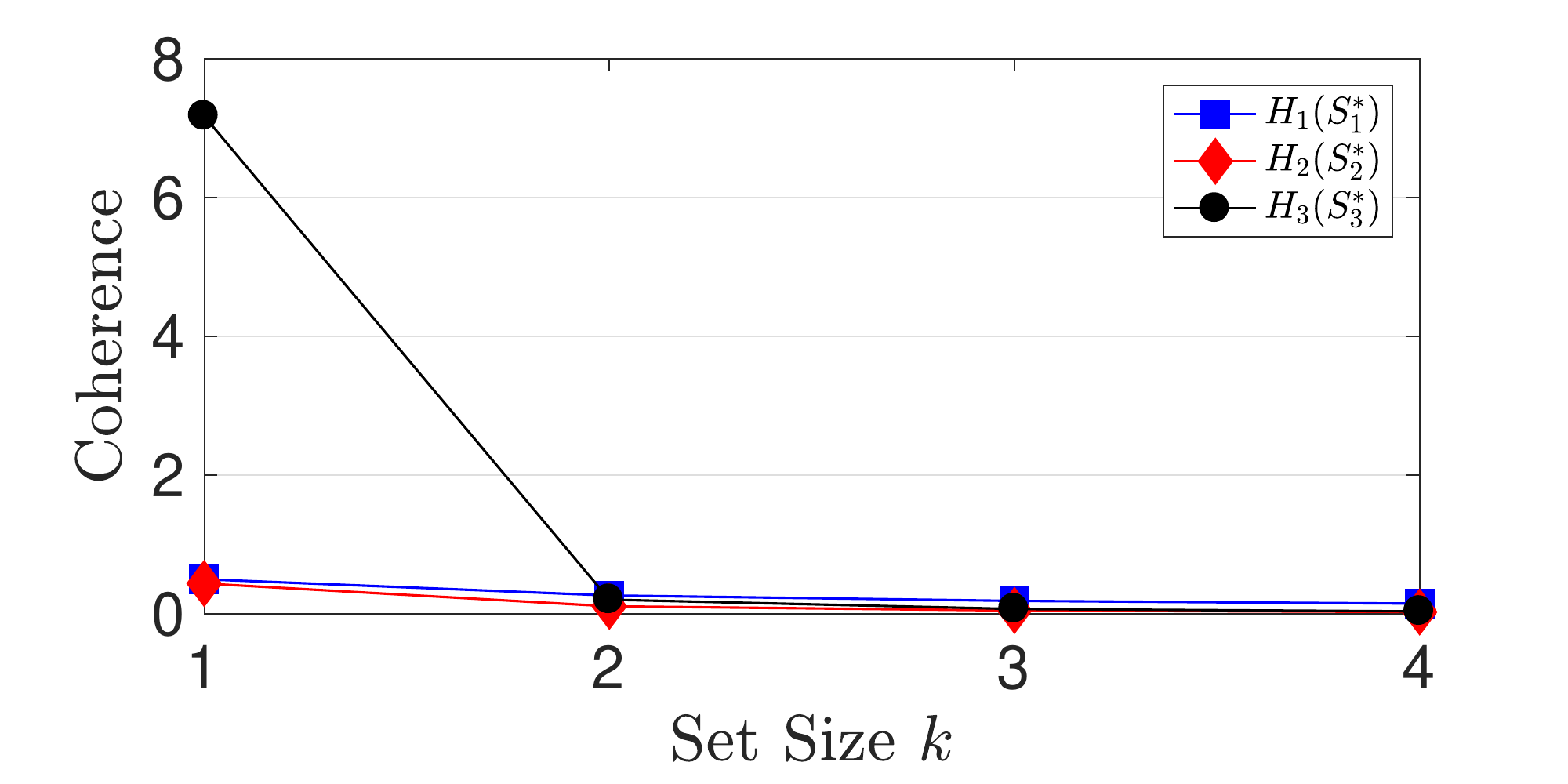}

\caption{Comparison of first, second, and third order optimal coherence in E-R graphs with $p=0.5$, $n=30$ over leader sets of size $k$.}
\label{O123OptCoh.fig}
\end{figure}

In Figure \ref{O123OptCoh.fig}, we compare the optimal coherence for first, second, and third order systems in Erd\H{o}s-R\'{e}nyi graphs of size $n=30$ with $p=0.5$, averaged over ten trials. We set $\Da=\I$ and each $a_i = a = \lceil \max{v \in V}  \frac{1}{\lambda_1(\Q{v})}\rceil$, so that ${a\lambda_i(\Q{S})>1}$ for $i=1,\ldots, n$. 

 As noted in (\ref{O23coh.eq}), $H_3(S) > H_2(S)$. Similarly, $H_1(S) > H_2(S)$, because $\sum_{i=1}^n\frac{1}{a\lambda_i(\Q{S})} >  \sum_{i=1}^n\frac{1}{(a\lambda_i(\Q{S}))^2}$. Since the optimal set minimizes the coherence, we can also note that the following is true
\begin{align*}
H_3(S_3^*) &> H_2(S_3^*) \geq H_2(S_2^*) \\
H_1(S_1^*) &> H_2(S_1^*) \geq H_2(S_2^*),
\end{align*}
as is demonstrated in the figure.

\begin{figure}
\centering

  \includegraphics[scale=.46]{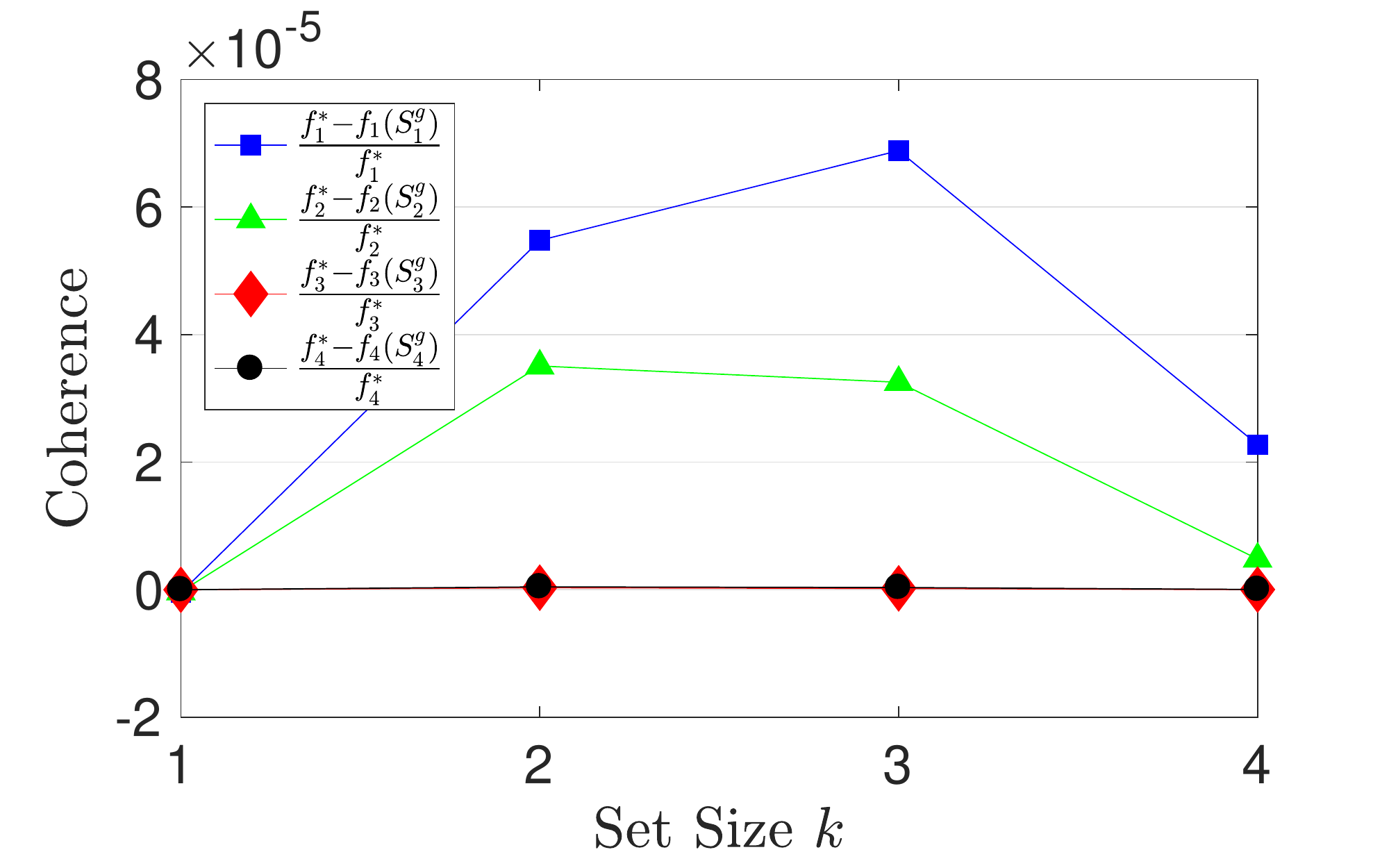}
  \centering
  \caption{Comparison of ratio from Theorem \ref{submodBound.thm} for an E-R graph with $p=0.5$, n=30.}

\label{cohRatio.fig}
\centering
\end{figure}

The submodularity of $f_m(S)$ guarantees that the coherence of a system of order ${m=1,2,3,4}$ with a set of leaders chosen by the greedy algorithm is within a constant bound of the optimal coherence.
Figure \ref{cohRatio.fig} plots the ratio $\frac{f_m^* -f_m(S_m^g)}{f_m^*}$ for $m=1,2,3,4$. This ratio has an upper bound of $\frac{1}{e}$ for each $m$. We compute this ratio for leader sets of size $k=1,2,3,4$ on Erd\H{o}s-R\'{e}nyi graphs with $p=0.5$ and $n=30$ and average the results over ten trials. We observe that for $f_1$, $f_2$, $f_3$, and $f_4$, not only is the ratio well below the upper bound guaranteed by Theorem \ref{cohBound.thm}, but it is in fact very close to $0$, and is $0$ when $k=1$, for all four systems. Thus, we can see that, for each $m$, the greedy algorithm generates a leader set that performs very close to optimal.

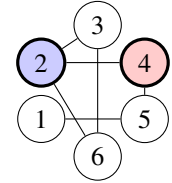
\begin{figure}
\centering
\begin{tikzpicture}[scale=0.5]

\node[shape=circle,draw=black] (5) at (1.25,-0.75) {5};
\node[shape=circle,draw=black, very thick, fill=blue!20] (2) at (-1.5,0.75) {2};
\node[shape=circle,draw=black] (3) at (0,1.75) {3};
\node[shape=circle,draw=black] (1) at (-1.5,-0.75) {1};
\node[shape=circle,draw=black, very thick, fill=red!20] (4) at (1.25,0.75) {4};
\node[shape=circle,draw=black] (6) at (0,-1.75) {6};

\foreach \from/\to in {1/5, 2/3, 2/4, 2/6, 3/6, 4/5}
    \draw (\from) -- (\to);
\end{tikzpicture}
\caption{Example graph where $\Da = I$ and all gains are set to $\lceil \max{v \in V}  \frac{1}{\lambda_1(\Q{v})}\rceil$. The optimal single leader for first order dynamics is node $2$, shaded in blue, while the optimal single leader for second and third order dynamics is $4$, shaded in red. }
\label{singleLeader.fig}
\end{figure}

Finally, we show that it is not guarateed that the same node is the optimal single leader for systems of different orders.
We only compare systems with $m=1$, $m=2$, and $m=3$ since it is possible to use the same gains $a_i$ in all three.
In Figure \ref{singleLeader.fig}, we show an example graph where $\Da=\I$ and each $a_i=a=\lceil \max{v \in V}  \frac{1}{\lambda_1(\Q{v})}\rceil$. In this case, selecting node $2$ as the leader minimizes the coherence under first order dynamics, but node $4$ is the best single leader under second and third order dynamics.

\section{Conclusion} \label{concl.sec}
We have studied leader-follower consensus dynamics in second, third, and fourth order systems.
First, we derived conditions under which such systems are stable, as well as proved that systems of order $m\geq 4$  are never stable when all gains are the same. 
We next developed expressions for coherence
in stable systems; these expressions are functions of the leader set.
We then proved that for $m=2,3,4$, a set function based on the coherence is submodular, and, therefore, the optimal set of leaders can be efficiently approximated using a greedy algorithm. Finally, we demonstrated the performance of higher order leader follower systems and the greedy algorithm through experiments and numerical examples.
In the future, we intend to study higher order consensus dynamics in directed graphs and graphs with negative edges.

\section{Appendix}
In this section, we present the proofs of Lemmas~\ref{O23.lem} and~\ref{O4.lem}.

\subsection{Proof of Lemma 1}

Our proof approach is based on a similar method to that used in \cite{SSLD15, MP17-2} and relies on the following theorem.
\begin{theorem}[\cite{L83}]
\label{submod}
A function $f:2^V \mapsto \mathbf{R}$ is submodular if and only if the derived set function $f_a:2^{V-\{a\}} \mapsto \mathbf{R}$, defined by
\begin{equation*}
f_a(X) = f(X \cup \{a\})-f(X),
\end{equation*}
is non-increasing for all $a \in V$.
\end{theorem}
\subsubsection{Proof that $f$ is non-decreasing}

Let $S_1 \subseteq S_2 \subseteq V$. First, we consider the case where ${S_1 = S_2 = \emptyset}$.  Then, ${f(S_1) = f(S_2) = 0}$, and so ${f(S_1) \leq f(S_2)}$, trivially.

Next, we consider the case where $S_1 = \emptyset$ and $S_2 \neq \emptyset$.
Since $f(S_1)=0$, we need only show that $f(S_2) \geq 0$.
By definition,
\begin{align*}
f(S_2) &= 2 \left(\max_{ s \in V}  \tr{\bQi{s}\bbQi{s}} \right)\\
&~~~~- \tr{\bQi{S_2}\bbQi{S_2}}\\
&\geq 2 \left( \max_{s \in S_2} \tr{\bQi{s}\bbQi{s}} \right)\\
&~~~~- \tr{\bQi{S_2}\bbQi{S_2}} \\
&=2 \sum_{i=1}^{n} \frac{1}{\lambda_i(b_1\Q{y})\lambda_i(b_2\Q{y}-b_3\I)}\\
&~~~~-\sum_{i=1}^{n} \frac{1}{\lambda_i(b_1\Q{S_2})\lambda_i(b_2\Q{S_2}-b_3\I)}\\
&=2 \sum_{i=1}^{n} \frac{1}{b_1\lambda_i(\Q{y})(b_2\lambda_i(\Q{y})-b_3)}\\
&~~~~-\sum_{i=1}^{n} \frac{1}{b_1\lambda_i(\Q{S_2})(b_2\lambda_i(\Q{S_2})-b_3)},
\end{align*}
where $y = \arg\max_{s \in S_2} \tr{\bQi{s}\bbQi{s}}$.  
Let $Z=S_2 \setminus \{y\}$. We write $\Q{S_2}$ as 
\begin{align*}
\Q{S_2} &= (\LL + \Da \D{S_2}) \\
&= \LL + \Da \D{Z} + \Da \D{y} \\
&= \Q{y}+\Da \D{Z}.
\end{align*}
The matrices $\Q{y}$ and $\Q{S_2}$ are positive definite \cite{RJME09}, and $\Da \D{Z}$ is a positive semidefinite matrix. Therefore, we can apply Weyl's Theorem, which gives us
\[ \lambda_i(\Q{S_2}) \geq \lambda_i(\Q{y})\]
and, thus, also
\[ \frac{1}{b_1\lambda_i(\Q{y})} \geq \frac{1}{b_1\lambda_i(\Q{S_2})}\]
and
\[\frac{1}{b_2\lambda_i(\Q{y})-b_3} \geq \frac{1}{b_2\lambda_i(\Q{S_2})-b_3}, \]
for $i=1,\ldots, n$. This implies that
\[\sum_{i=1}^{n} \frac{1}{b_1\lambda_i(\Q{S_2})(b_2\lambda_i(\Q{S_2})-b_3)} \leq \sum_{i=1}^{n} \frac{1}{b_1\lambda_i(\Q{y})(b_2\lambda_i(\Q{y})-b_3)}\]
and, therefore, $f(S_2) \geq 0$.

Finally, we consider the case where $S_1 \neq \emptyset$ and $S_2 \neq \emptyset$. Then,
\begin{align*}
f(S_1) - f(S_2) &=\tr{\bQi{S_2}\bbQi{S_2}} \\
&~~~~- \tr{\bQi{S_1}\bbQi{S_1}} \\
&=\sum_{i=1}^{n} \frac{1}{b_1\lambda_i(\Q{S_2})(b_2\lambda_i(\Q{S_2})-b_3)}\\
&~~~~-\sum_{i=1}^{n} \frac{1}{b_1\lambda_i(\Q{S_1})(b_2\lambda_i(\Q{S_1})-b_3)},
\end{align*}

Let $Z = S_2 \setminus S_1$ so that $\Q{S_2} = \Q{S_1} + \Da \D{Z}$.
By the same argument as in the previous case, we can use Weyl's Theorem to find that
\[ \lambda_i(\Q{S_2}) \geq \lambda_i(\Q{S_1}), \]
for $i=1,\ldots, n$.
We can then conclude that 
\begin{align*}
&\sum_{i=1}^{n} \frac{1}{b_1\lambda_i(\Q{S_2})(b_2\lambda_i(\Q{S_2})-b_3)} \\
&~~~~\leq \sum_{i=1}^{n} \frac{1}{b_1\lambda_i(\Q{S_1})(b_2\lambda_i(\Q{S_2})-b_3)}
\end{align*}
and, therefore, $f(S_1) \leq f(S_2)$.
Thus, $f$ is non-decreasing.

\subsection{Proof that $f$ is submodular}
To prove that $f$ is submodular, we first define the set function ${f_a : 2^{V \setminus \{a\} } \mapsto \mathbf{R}}$,
\begin{equation*}
f_a(S) =  f( S \cup \{a\}) - f(S), 
\end{equation*}
and show that it is non-increasing.

Let $S_1 \subseteq S_2 \subseteq V \setminus \{a\}$. 
First, we consider the case where $S_1=S_2=\emptyset$.
In this case, $f_a(S_1) = f_a(S_2) = f(\{a\})$, so, trivially, $f_a(S_1) \geq f_a(S_2)$ and $f_{a}$ is non-increasing for $S = \emptyset$. \\

Next, we consider the case where $S_1 = \emptyset$, $S_2 \neq \emptyset$.
Then,
\begin{align*}
&f_a(S_1) - f_a(S_2) = f(\{a\}) - \big(f(S_2 \cup \{a\}) - f(S_2) \big) \\
&= \left(C - \tr{\bQi{a}\bbQi{a}}\right)  \\
&~~~~- \tr{\bQi{S_2}\bbQi{S_2}} \\
&~~~~+ \tr{\bQi{S_2 \cup \{a\}}\bbQi{S_2 \cup \{a\}}}. 
\end{align*}

Recall that $C = 2 \left(\max_{ s \in V}  \tr{\bQi{s}\bbQi{S}}\right)$, and so
\begin{align*}
 &C - \tr{\bQi{a}\bbQi{a}} \\
 &\geq \max_{s \in V}  \tr{\bQi{s}\bbQi{s}} \\
 &\geq \max_{s \in S_2} \tr{\bQi{s}\bbQi{s}}.
\end{align*}

We have already shown that $f$ is non-decreasing. Thus,
\begin{align*}
&\max_{s\in S_2} \tr{\bQi{s}\bbQi{s}} \\
&~~~~- \tr{\bQi{S_2}\bbQi{S_2}} \geq 0.
\end{align*}
Therefore, $f_a(S_1) - f_a(S_2) \geq 0$.\\

Finally, we consider the case where $S_1 \neq \emptyset$ and $S_2 \neq \emptyset$.
Let $\D{a}$ be a diagonal matrix with $\D{a}(a,a) = 1$ and all other entries equal to 0. Then $f_a(S)$ becomes:
\begin{align*}
&f_a(S) = C - \tr{\bQi{S \cup \{a\}}\bbQi{S \cup \{a\}}} \\
&~~~~- \left(C -  \tr{\bQi{S}\bbQi{S}}\right)\\
&~~=- \tr{\big(b_1(\Q{S} + \Da \D{a})\big)^{-1}\big(b_2(\Q{S} +\Da \D{a})-b_3\I\big)^{-1}}\\
&~~~~+ \tr{\bQi{S}\bbQi{S}}.
\end{align*}

We define the functions $\F (t)$ over $t \in [0,1]$, as:
\begin{align*}
 \F (t) &= \Q{S_1} + t(\Q{S_2}-\Q{S_1}).
 \end{align*}
 Note that $\F(0) = \Q{S_1}$ and $\F(1) = \Q{S_2}$. Let
\begin{align}
&\hat{f}_{a} (\F (t)) =\nonumber \\
& -\tr{\big(b_1(\F (t)+ \Da \D{a})\big)^{-1}\big(b_2(\F (t)+\Da \D{a}\big) -b_3\I)^{-1}} \nonumber  \\
&~~~~+\tr{\big(b_1\F (t)\big)^{-1}\big(b_2\F (t)-b_3\I\big)^{-1}}. \nonumber 
\end{align}

We then take the derivative of $\hat{f}_a$ with respect to $t$,
\begin{align}
&\frac{d}{dt}\hat{f}_{a}(\F (t))=  \frac{d}{dt}  \Big( -\textbf{tr}\Big(\big(b_1(\F (t)+ \Da \D{a})\big)^{-1}\nonumber \\
&~~~~~~\big(b_2\left(\F (t)+\Da \D{a}\right) -b_3I\big)^{-1}\Big) \nonumber  \\
&~~~~+\tr{\big(b_1\F (t)\big)^{-1}\big(b_2\F (t)-b_3\I\big)^{-1}} \Big). \label{deriv1}
\end{align}

Note that
\begin{align} 
\scalemath{0.95}{\frac{d}{dt}(\A^{-1}\B^{-1}) = -\A^{-1}\frac{d}{dt}(\A)\A^{-1}\B^{-1}-\A^{-1}\B^{-1}\frac{d}{dt}(\B)\B^{-1}.} \label{matDeriv.eq}
\end{align}
By applying (\ref{matDeriv.eq}) and the cyclic property of the trace to (\ref{deriv1}), we obtain:
%
\begin{align}
&\frac{d}{dt}\hat{f}_{a}(\F (t)) \label{deriv2}\\
&=\textbf{tr}\Big(\big(b_1(\F (t)+ \Da \D{a})\big)^{-1}\big(b_2\left(\F (t)+\Da \D{a}\right) -b_3\I\big)^{-1} \nonumber \\
&~~~~\big(b_1(\F (t)+ \Da \D{a})\big)^{-1}\big(b_1(\Q{S_2}-\Q{S_1})\big)\Big) \label{deriv2a} \\
 &~~+ \textbf{tr}\Big(\big(b_2\left(\F (t)+\Da \D{a}\right) -b_3\I\big)^{-1}\big(b_1(\F (t)+ \Da \D{a})\big)^{-1}\nonumber \\
 &~~~~\big(b_2\left(\F (t)+\Da \D{a}\right) -b_3\I\big)^{-1}\big(b_2(\Q{S_2}-\Q{S_1})\big)\Big) \label{deriv2b} \\
 &~~-\textbf{tr}\Big(\big(b_1\F (t)\big)^{-1}\big(b_2\F (t)-b_3\I\big)^{-1}\big(b_1\F (t)\big)^{-1}\nonumber \\
 &~~~~~~~~\big(b_1(\Q{S_2}-\Q{S_1})\big)\Big) -\textbf{tr}\Big(\big(b_2\F (t)-b_3\I\big)^{-1}\nonumber \\
 &~~~~~~~~\big(b_1\F (t)\big)^{-1}\big(b_2\F (t)-b_3\I\big)^{-1}\big(b_2(\Q{S_2}-\Q{S_1})\big)\Big). \nonumber 
\end{align}
We now show that (\ref{deriv2}) is non-positive. To do so, we first must simplify the expression.
We define $\daa$ as the vector of all zeros except the $a^{th}$ component, which has value $\sqrt{\vaa}$. We also define
\[\scalemath{0.95}{\W=\frac{\big(b_1\F (t)\big)^{-1}b_1\Da \D{a} \big(b_1\F (t)\big)^{-1}}{1 - b_1\daa^T \big(b_1\F (t)\big)^{-1} \daa}}\]
and
\[\scalemath{0.95}{\Z=\frac{\big(b_2\F (t)-b_3\I\big)^{-1}b_2\Da \D{a} \big(b_2\F (t)-b_3\I\big)^{-1}}{1 - b_2\daa^T \big(b_2\F (t) -b_3\I\big)^{-1}\daa}}.\]

Using $\W$ and $\Z$, we apply the Sherman-Morrison formula to the first factor of (\ref{deriv2a}) to obtain:

\begin{align}
&\Big(\big(b_1\F (t)\big)^{-1}-\W\Big)\Big(\big(b_2\F (t)-b_3\I\big)^{-1}-\Z\Big)\Big(\big(b_1\F (t)\big)^{-1}-\W\Big) \nonumber \\
&=\big(b_1\F (t)\big)^{-1}\big(b_2\F (t)-b_3\I\big)^{-1}\big(b_1\F (t)\big)^{-1} \nonumber \\
&~~~~- \big(b_1\F (t)\big)^{-1}\Big(\big(b_2\F (t)-b_3\I\big)^{-1}-\Z\Big)\W \nonumber \\
&~~~~ -\Big(\big(b_1\F (t)\big)^{-1}-\W\Big)\Z\big(b_1\F (t)\big)^{-1}  \nonumber \\
&~~~~- \W \big(b_2\F (t)-b_3\I\big)^{-1}\Big(\big(b_1\F (t)\big)^{-1}-\W\Big)- \W\Z\W\nonumber \\
&=\big(b_1\F (t)\big)^{-1}\big(b_2\F (t)-b_3\I\big)^{-1}\big(b_1\F (t)\big)^{-1} \nonumber \\
&~~~~- \big(b_1\F (t)\big)^{-1}\big(b_2\left(\F (t)+\Da \D{a}\right) -b_3\I\big)^{-1}\W \nonumber \\
&~~~~ -\big(b_1(\F (t)+ \Da \D{a})\big)^{-1}\Z\big(b_1\F (t)\big)^{-1}  \nonumber \\
&~~~~- \W \big(b_2\F (t)-b_3\I\big)^{-1}\big(b_1(\F (t)+ \Da \D{a})\big)^{-1}- \W\Z\W. \label{deriv2afactor}
\end{align}

We apply the Sherman-Morrison formula to the first factor of (\ref{deriv2b}) as well to obtain:
\begin{align}
&\Big(\big(b_2\F (t)-b_3\I\big)^{-1}-\Z\Big)\Big(\big(b_1\F (t)\big)^{-1}-\W\Big)\nonumber \\
&~~~~\Big(\big(b_2\F (t)-b_3\I\big)^{-1}-\Z\Big)\nonumber \\
&=\big(b_2\F (t)-b_3\I\big)^{-1}\big(b_1\F (t)\big)^{-1}\big(b_2\F (t)-b_3\I\big)^{-1} \nonumber \\
&~~~~- \big(b_2\F (t)-b_3\I\big)^{-1}\big(b_1(\F (t)+ \Da \D{a})\big)^{-1}\Z\nonumber \\
&~~~~- \big(b_2\left(\F (t)+\Da \D{a}\right) -b_3\I\big)^{-1}\W\big(b_2\F (t)-b_3\I\big)^{-1}\nonumber \\
& ~~~~- \Z\big(b_1\F (t)\big)^{-1}\big(b_2\left(\F (t)+\Da \D{a}\right) -b_3\I\big)^{-1}-\Z\W\Z. \label{deriv2bfactor}
\end{align}

We can then substitute (\ref{deriv2afactor}) and (\ref{deriv2bfactor}) back into (\ref{deriv2}) to get
\begin{align}
&\textbf{tr}\Bigg(\Big(-\big(b_1\F (t)\big)^{-1}\big(b_2\left(\F (t)+\Da \D{a}\right) -b_3\I\big)^{-1}\W \nonumber \\
& -\big(b_1(\F (t)+ \Da \D{a})\big)^{-1}\Z\big(b_1\F (t)\big)^{-1}  \nonumber \\
&- \W \big(b_2\F (t)-b_3\I\big)^{-1}\big(b_1(\F (t)+ \Da \D{a})\big)^{-1}- \W\Z\W\Big) \nonumber \\
&~~~~\big(b_1(\Q{S_2}-\Q{S_1})\big)\Bigg) \nonumber\\
 &~~+ \textbf{tr}\Bigg( \Big(- \big(b_2\F (t)-b_3\I\big)^{-1}\big(b_1(\F (t)+ \Da \D{a})\big)^{-1}\Z\nonumber \\
&- \big(b_2\left(\F (t)+\Da \D{a}\right) -b_3\I\big)^{-1}\W\big(b_2\F (t)-b_3\I\big)^{-1}\nonumber \\
& - \Z\big(b_1\F (t)\big)^{-1}\big(b_2\left(\F (t)+\Da \D{a}\right) -b_3\I\big)^{-1}-\Z\W\Z\Big)\nonumber \\
&~~~~\big(b_2(\Q{S_2}-\Q{S_1})\big)\Bigg). \label{deriv3}
\end{align}

Note that $b_1\F (t)$, $b_1(\F (t) + \Da \D{a})$ $b_2\F (t)-b_3\I$, and $b_2(\F (t) + \Da \D{a})-b_3\I$ are all grounded Laplacians, and therefore their inverses are element-wise positive for all possible values of $t$ and $a$~\cite{M93}.

We also note that $\W$ and $\Z$ are both element-wise non-negative: the numerator is the product of two element-wise non-negative matrices  and two element-wise positive matrices and the denominator is the positive scalar $1 - \vaa \F (t)^{-1}_{(a,a)}$, $1 - b_2 \vaa (b_2\F (t)-b_3\I)^{-1}_{(a,a)}$, respectively. We can thus see that every element of the first factor in each trace in (\ref{deriv3}) is the product of three element-wise non-negative matrices.  Therefore, the first factor in both traces is element-wise non-positive.
Finally, we note that 
\begin{align*}
\Q{S_2}-\Q{S_1} &= (\LL + \Da \D{S_2}) - (\LL + \Da \D{S_1})  \\
&~~~~~=\Da (\D{S_2} - \D{S_1}). 
\end{align*}
Thus, $\Q{S_2}-\Q{S_1}$ is a diagonal matrix where the $(i,i)^{th}$ component is positive if ${i \in S_2 \setminus S_1}$ and $0$ otherwise.
Therefore, (\ref{deriv3}), as the sum of two traces, each the product of an element-wise non-positive matrix  and an element-wise non-negative matrix, is element-wise non-positive. Thus, $\frac{d}{dt} \hat{f}_a(\F (t)) \leq 0$.
 
Now, consider the following equality:
\begin{align*}
&\hat{f}_a(F(1)) = \\
&~~~~\hat{f}_a(F(0))+\int_0^1  \textstyle \frac{d}{dt}\hat{f}_a(\F (t))dt.
\end{align*}
Note that $\hat{f}_a(F(1)) = f_{a}(S_2)$ and $\hat{f}_a(F(0)) = f_{a}(S_1)$, and, as shown above,  $\frac{d}{dt}\hat{f_a}(\F (t))dt$ is non-positive. 
Thus,  $ f_{a}(S_1) \geq f_{a}(S_2)$. 

 From this we can conclude that $f_{a}$ is non-increasing for all $a \in V$, and, therefore, by Theorem~\ref{submod}, $f$ is submodular.

\subsection{Proof of Lemma 2}
Lemma~\ref{O4.lem} is proved using similar methods to those used in the proof of Lemma \ref{O23.lem}.

\subsubsection{Proof that $f$ is non-decreasing}
Let $S_1 \subseteq S_2 \subseteq V$. First, we consider the case where $S_1 = S_2 = \emptyset$.  Then, $f(S_1) = f(S_2) = 0$, and so $f(S_1) \leq f(S_2)$, trivially.

Next, we consider the case where $S_1 = \emptyset$ and $S_2 \neq \emptyset$.
Since $f(S_1)=0$, we need only show that $f(S_2) \geq 0$.
By definition,
\begin{align*}
f(S_2) &= 2 \left(\max_{ s \in V}  \tr{\Qit{s}(b_1\Q{s}-\I)\big((b_1-b_2)\Q{s}-\I\big)^{-1}} \right)\\
&~~~~- \tr{\Qit{S_2}(b_1\Q{S_2}-\I)\big((b_1-b_2)\Q{S_2}-\I\big)^{-1}}\\
&\geq 2 \left( \max_{s \in S_2} \tr{\Qit{s}(b_1\Q{s}-\I)\big((b_1-b_2)\Q{s}-\I\big)^{-1}} \right)\\
&~~~~- \tr{\Qit{S_2}(b_1\Q{S_2}-\I)\big((b_1-b_2)\Q{S_2}-\I\big)^{-1}} \\
&=2 \sum_{i=1}^{n} \frac{b_1\lambda_i(\Q{y})-1}{\lambda_i(\Q{y})^2((b_1-b_2)\lambda_i(\Q{y})-1)}\\
&~~~~-\sum_{i=1}^{n} \frac{b_1\lambda_i(\Q{S_2})-1}{\lambda_i(\Q{S_2})^2((b_1-b_2)\lambda_i(\Q{S_2})-1)},
\end{align*}
where $y = \arg\max_{s \in S_2} \tr{\Qit{s}(b_1\Q{s}-\I)\big((b_1-b_2)\Q{s}-\I\big)^{-1}}$. 

Let $Z=S_2 \setminus \{y\}$. We write $\Q{S_2}$ as 
\begin{align*}
\Q{S_2} &= (\LL + \Da \D{S_2}) \\
&= \LL + \Da \D{Z} + \Da \D{y} \\
&= \Q{y}+\Da \D{Z}.
\end{align*}
The matrices $\Q{y}$ and $\Q{S_2}$ are positive definite \cite{RJME09}, and $\Da \D{Z}$ is a positive semidefinite matrix. Therefore, we can apply Weyl's Theorem, which gives us
\[ \lambda_i(\Q{y}) \leq \lambda_i(\Q{S_2}),\]
and thus also 
\begin{align}
\frac{1}{\lambda_i(\Q{y})} \geq \frac{1}{\lambda_1(\Q{S_2})}. \label{step1.eq}
\end{align}
To show that $f(S_2)\geq 0$, we need to also show that 
\begin{align*}
\frac{b_1\lambda_i(\Q{y})-1}{(b_1-b_2)\lambda_i(\Q{y})-1} &\geq \frac{b_1\lambda_i(\Q{S_2})-1}{(b_1-b_2)\lambda_i(\Q{S_2})-1}
\end{align*}
holds for all $i=1, \ldots, n$. To show this, we first multiply both sides of (\ref{step1.eq}) by $\frac{1}{b_1}$ and note that 
\begin{align*}
\frac{1}{b_1\lambda_i(\Q{y})} &\geq \frac{1}{b_1\lambda_i(\Q{S_2})} 
\end{align*}
is true if and only if 
\begin{align}
\frac{1}{1-\frac{1}{b_1\lambda_i(\Q{y})}} &\geq \frac{1}{1- \frac{1}{b_1\lambda_1(\Q{S_2})}}. \label{step3.eq}
\end{align}
We rearrange (\ref{step3.eq}) and multiply both sides by $\frac{b_2}{b_1}$ to get:
\begin{align}
\frac{b_2\lambda_i(\Q{y})}{b_1\lambda_i(\Q{y})-1} &\geq \frac{b_2\lambda_i(\Q{S_2})}{b_1\lambda_i(\Q{S_2})- 1}, \nonumber 
\end{align}
which we again note is true if and only if 
\begin{align*}
\frac{1}{1-\frac{b_2\lambda_i(\Q{y})}{b_1\lambda_i(\Q{y})-1}} &\geq \frac{1}{1- \frac{b_2\lambda_i(\Q{S_2})}{b_1\lambda_i(\Q{S_2})- 1}} 
\end{align*}
is true.
We then rearrange the expression to get
\begin{align*}
\frac{b_1\lambda_i(\Q{y})-1}{(b_1-b_2)\lambda_i(\Q{y})-1} &\geq \frac{b_1\lambda_i(\Q{S_2})-1}{(b_1-b_2)\lambda_i(\Q{S_2})-1}.
\end{align*}

From this, we can conclude that 
\begin{align*}
&2 \sum_{i=1}^{n} \frac{b_1\lambda_i(\Q{y})-1)}{\lambda_i(\Q{y})^2((b_1-b_2)\lambda_i(\Q{y})-1)}\\
&~~~~\geq \sum_{i=1}^{n} \frac{b_1\lambda_i(\Q{S_2})-1)}{\lambda_i(\Q{S_2})^2((b_1-b_2)\lambda_i(\Q{S_2})-1)}
\end{align*}
and, therefore, $f(S_2) \geq 0$.

Finally, we consider the case where $S_1 \neq \emptyset$ and $S_2 \neq \emptyset$. Then,
\begin{align*}
&f(S_1) - f(S_2) \\
&~~=\tr{\Qit{S_2}(b_1\Q{S_2}-\I)\big((b_1-b_2)\Q{S_2}-\I\big)^{-1}} \\
&~~~~~~- \tr{\Qit{S_1}(b_1\Q{S_1}-\I)\big((b_1-b_2)\Q{S_1}-\I\big)^{-1}}\\
&~~=\sum_{i=1}^{n} \frac{b_1\lambda_i(\Q{S_2})-1)}{\lambda_i(\Q{S_2})^2((b_1-b_2)\lambda_i(\Q{S_2})-1)}\\
&~~~~~~-\sum_{i=1}^{n} \frac{b_1\lambda_i(\Q{S_1})-1)}{\lambda_i(\Q{S_1})^2((b_1-b_2)\lambda_i(\Q{S_1})-1)}.
\end{align*}

Let $Z = S_2 \setminus S_1$ so that $\Q{S_2} = \Q{S_1} + \Da \D{Z}$.
By the same arguments as used in the previous case, we find that
\[ \frac{1}{\lambda_i(\Q{S_1})} \geq \frac{1}{\lambda_i(\Q{S_2})}, \]
and 
\begin{align*}
\frac{b_1\lambda(\Q{S_1})-1}{(b_1-b_2)\lambda(\Q{S_1})-1} &\geq \frac{b_1\lambda(\Q{S_2})-1}{(b_1-b_2)\lambda(\Q{S_2})-1}
\end{align*}
for $i=1, \ldots, n$.
Thus, 
\begin{align*} 
&\sum_{i=1}^{n} \frac{b_1\lambda_i(\Q{S_1})-1)}{\lambda_i(\Q{S_1})^2((b_1-b_2)\lambda_i(\Q{S_1)-1)}}\\
&~~~~\geq \sum_{i=1}^{n} \frac{b_1\lambda_i(\Q{S_2})-1)}{\lambda_i(\Q{S_2})^2((b_1-b_2)\lambda_i(\Q{S_2})-1)}
\end{align*}
and $f(S_1) \leq f(S_2)$.
Therefore, $f$ is non-decreasing.

\subsubsection{Proof that $f$ is submodular}

To prove that $f$ is submodular, we  define the set function ${f_a : 2^{V \setminus \{a\} } \mapsto \mathbf{R}}$,
\begin{equation*}
f_a(S) =  f( S \cup \{a\}) - f(S), 
\end{equation*}
and show that it is non-increasing.

Let $S_1 \subseteq S_2 \subseteq V \setminus \{a\}$. 
First, we consider the case where $S_1=S_2=\emptyset$.
In this case, $f_a(S_1) = f_a(S_2) = f(\{a\})$, so, trivially, $f_a(S_1) \geq f_a(S_2)$ and $f_{a}$ is non-increasing for $S = \emptyset$. \\

Next, we consider the case where $S_1 = \emptyset$, $S_2 \neq \emptyset$.
Then,
\begin{align*}
&f_a(S_1) - f_a(S_2) = f(\{a\}) - \big(f(S_2 \cup \{a\}) - f(S_2) \big) \\
&= \left(C - \tr{\Qit{a}(b_1\Q{a}-\I)\big((b_1-b_2)\Q{a}-\I\big)^{-1}}\right)  \\
&~~~~- \tr{\Qit{S_2}(b_1\Q{S_2}-\I)\big((b_1-b_2)\Q{S_2}-\I\big)^{-1}} \\
&~~~~+ \tr{\Qit{S_2\cup \{a\}}(b_1\Q{S_2 \cup \{a\}}-\I)\big((b_1-b_2)\Q{S_2 \cup \{a\}}-\I\big)^{-1}}. 
\end{align*}

Recall that \[C = 2 \left(\max_{ s \in V}  \tr{\Qit{s}(b_1\Q{s}-\I)\big((b_1-b_2)\Q{s}-\I\big)^{-1}}\right),\] and so
\begin{align*}
 &C - \tr{\Qit{a}(b_1\Q{a}-\I)\big((b_1-b_2)\Q{a}-\I\big)^{-1}} \\
 &\geq \max_{s \in V}  \tr{\Qit{s}(b_1\Q{s}-\I)\big((b_1-b_2)\Q{s}-\I\big)^{-1}} \\
 &\geq \max_{s \in S_2} \tr{\Qit{s}(b_1\Q{s}-\I)\big((b_1-b_2)\Q{s}-\I\big)^{-1}}.
\end{align*}

We have already shown that $f$ is non-decreasing. Thus, 
\begin{align*}
&\max_{s\in S_2} \tr{\Qit{s}(b_1\Q{s}-\I)\big((b_1-b_2)\Q{s}-\I\big)^{-1}} \\
&~~~~- \tr{\Qit{S_2}(b_1\Q{S_2}-\I)\big((b_1-b_2)\Q{S_2}-\I\big)^{-1}} \geq 0.
\end{align*}
Therefore, $f_a(S_1) - f_a(S_2) \geq 0$.\\

Finally, we consider the case where $S_1 \neq \emptyset$ and $S_2 \neq \emptyset$.
Let $\D{a}$ be a diagonal matrix with $\D{a}(a,a) = 1$ and all other entries equal to 0. Then $f_a(S)$ becomes:
\begin{align*}
&f_a(S)\\
& ~~=C - \tr{\Qit{S \cup \{a\}}(b_1\Q{S\cup \{a\}}-\I)\big((b_1-b_2)\Q{S \cup \{a\}}-\I\big)^{-1}}\\
&~~~~- \left(C -  \tr{\Qit{S}(b_1\Q{S}-\I)\big((b_1-b_2)\Q{S}-\I\big)^{-1}}\right)\\
&~~=- \textbf{tr}\Big(\left((\Q{S} + \Da \D{a})\right)^{-2}\big(b_1(\Q{S} +\Da \D{a})-\I\big)\\
&~~~~~~~~\big((b_1-b_2)(\Q{S}+\Da \D{a})-\I\big)^{-1}\Big)\\
&~~~~+ \tr{\Qit{S}(b_1\Q{S}-\I)\big((b_1-b_2)\Q{S}-\I\big)^{-1}}.
\end{align*}

We define the functions $\F (t)$ over $t \in [0,1]$, as:
\begin{align*}
 \F (t) &= \Q{S_1} + t(\Q{S_2}-\Q{S_1}).
 \end{align*}
 Note that $F(0) = \Q{S_1}$ and $F(1) = \Q{S_2}$. Let
\begin{align}
&\hat{f}_{a} (\F (t)) = \nonumber \\
& -\textbf{tr}\Big(\left(\F (t)+ \Da \D{a}\right)^{-2}\big(b_1(\F (t)+\Da \D{a}) -\I\big)\nonumber \\
&~~~~~~~~\big((b_1-b_2)(\F (t)+\Da \D{a})-\I\big)^{-1}\Big) \nonumber  \\
&~~~~+\tr{\F (t)^{-2}(b_1\F (t)-\I)\big((b_1-b_2)\F (t)-\I\big)^{-1}} \nonumber \\
& =-\textbf{tr}\Big(\big(\F (t)+ \Da \D{a}\big)^{-2} + b_2\big(\F (t)+\Da \D{a}\big)^{-1}\nonumber \\
&~~~~~~\big((b_1-b_2)(\F (t)+\Da \D{a}) -\I\big)^{-1}\Big)\nonumber \\
&~~~~+\tr{\F (t)^{-2}+b_2\F (t)^{-1}\big((b_1-b_2)\F (t)-\I\big)^{-1}}. \label{fahat4.eq}
\end{align}

We use the formula for the derivative of  $\A^{-1}\B^{-1}$ (\ref{matDeriv.eq}) and the cyclic property of the trace to take the derivative of (\ref{fahat4.eq}) with respect to $t$.

\begin{align}
&\frac{d}{dt}\hat{f}_{a}(\F (t)) \label{deriv4}\\
&=\textbf{tr}\Bigg( (\Q{S_2}-\Q{S_1})\Big[2\big(\F (t)+ \Da \D{a}\big)^{-3} \label{deriv41} \\
&~~+\Big(b_2\big(\F (t)+ \Da \D{a}\big)^{-1}\label{deriv42} \\
&~~~~~~\big[(b_1-b_2)(\F (t)+\Da \D{a}) -\I\big]^{-1}\big(\F (t)+ \Da \D{a}\big)^{-1}\Big) \nonumber \\
&~~+\Big(b_2(b_1-b_2)\big[(b_1-b_2)(\F (t)+\Da \D{a}) -\I\big]^{-1}  \label{deriv43} \\
&~~~~~~(\F (t)+ \Da \D{a})^{-1}\big[(b_1-b_2)(\F (t)+\Da \D{a}) -\I\big]^{-1} \Big)  \nonumber\\
&~~-2\F (t)^{-3} -\Big(b_2\F (t)^{-1}\big[(b_1-b_2)\F (t) -\I\big]^{-1}\F (t)^{-1}\Big) \nonumber \\
&~~~~-\Big(b_2(b_1-b_2)\big[(b_1-b_2)\F (t) -\I\big]^{-1}\F (t)^{-1}\nonumber \\
&~~~~~~\big[(b_1-b_2)\F (t) -\I\big]^{-1} \Big) \Big] \Bigg). \nonumber
\end{align}
We now show that (\ref{deriv4}) is non-positive. To do so, we first simplify the expression by expanding components (\ref{deriv41}), (\ref{deriv42}), and (\ref{deriv43}).

We define $\daa$ as the vector of all zeros except the $a^{th}$ component, which has value $\sqrt{\vaa}$. For simplicity, we define
\[ \G(t)=(b_1-b_2)\F (t) -\I.\]
We also define
\[\W=\frac{\F (t)^{-1}\Da \D{a} \F (t)^{-1}}{1 - \daa^T \F (t)^{-1} \daa}\]
and
\begin{align*}
&\Z=\\
&\frac{\big((b_1-b_2)\F (t)-\I\big)^{-1}(b_1-b_2)\Da \D{a} \big((b_1-b_2)\F (t)-\I\big)^{-1}}{1 - (b_1-b_2)\daa^T \big((b_1-b_2)\F (t) -\I\big)^{-1}\daa}.
\end{align*}

Using $\W$ and $\Z$, we apply the Sherman-Morrison formula to (\ref{deriv41}) to obtain:

\begin{align}
&(\F (t)+ \Da \D{a})^{-3}=(\F (t)^{-1} -\W)^{3} \nonumber \\
&= \F (t)^{-3} -\F (t)^{-1} (\F (t)^{-1}-\W)\W \nonumber \\
&~~~~-(\F (t)^{-1}-\W)\W \F (t)^{-1}\nonumber \\
&~~~~-\W(\F (t)^{-1}-\W)\F (t)^{-1}-\W^{3}\nonumber \\
&=\F (t)^{-3} -\F (t)^{-1} (\F (t)+ \Da \D{a})^{-1}\W\nonumber \\
&~~~~-(\F (t)+ \Da \D{a})^{-1}\W \F (t)^{-1}\nonumber \\
&~~~~-\W(\F (t)+ \Da \D{a})^{-1}\F (t)^{-1}\nonumber \\
&~~~~-\W^{3}.  \label{deriv41factor}
\end{align}

We then simplify (\ref{deriv42}) and (\ref{deriv43}) in the same way to get:
\begin{align}
&(\F (t)+ \Da \D{a})^{-1}\big[(b_1-b_2)(\F (t)+\Da \D{a}) -\I\big]^{-1}\nonumber \\
&~~~~~~(\F (t)+ \Da \D{a})^{-1}\Big) \nonumber \\
&=(\F (t)^{-1}-\W)(\G (t)^{-1}-\Z)(\F (t)^{-1}-\Z) \nonumber \\
&=\F (t)^{-1}\G (t)^{-1}\F (t)^{-1}\nonumber \\
&~~~~-\F (t)^{-1}\big[(b_1-b_2)(\F (t)+\Da \D{a}) -\I\big]^{-1}\W\nonumber \\
&~~~~-(\F (t)+ \Da \D{a})^{-1}\Z \F (t)^{-1}-\W \G (t)^{-1}(\F (t)+ \Da \D{a})^{-1}\nonumber \\
&~~~~-\W\Z\W. \label{deriv43factor}
\end{align}

and

\begin{align}
&\big[(b_1-b_2)(\F (t)+\Da \D{a}) -\I\big]^{-1}(\F (t)+ \Da \D{a})^{-1}\nonumber \\
&~~~~\big[(b_1-b_2)(\F (t)+\Da \D{a}) -\I\big]^{-1} \nonumber \\
&= (\G (t)^{-1}-\Z)(\F (t)^{-1}-\W)(\G (t)^{-1}-\Z)\nonumber \\
&=\G (t)^{-1}\F (t)^{-1}\G (t)^{-1} - \G (t)^{-1}(\F (t)+ \Da \D{a})^{-1}\Z \nonumber \\
&~~~~- \big[(b_1-b_2)(\F (t)+\Da \D{a}) -\I\big]^{-1}\W \G (t)^{-1}\nonumber \\
&~~~~-\Z \F (t)^{-1}\big[(b_1-b_2)(\F (t)+\Da \D{a}) -\I\big]^{-1}\nonumber \\
&~~~~-\Z\W\Z.  \label{deriv42factor}
\end{align}

We can then substitute (\ref{deriv41factor}), (\ref{deriv42factor}), and (\ref{deriv43factor}) back into (\ref{deriv2}) to get
\begin{align}
&\frac{d}{dt} \hat{f}_a(\F (t))\nonumber \\
&~~=\textbf{tr}\Bigg( (\Q{S_2}-\Q{S_1})\Big[2\Big(-\F (t)^{-1} (\F (t)+ \Da \D{a})^{-1}\W\nonumber \\
&~~~~~~~~-(\F (t)+ \Da \D{a})^{-1}\W \F (t)^{-1}\nonumber \\
&~~~~~~~~-\W(\F (t)+ \Da \D{a})^{-1}\F (t)^{-1}-\W^{3}\Big) \nonumber\\
&~~~~+ b_2\Big(-\F (t)^{-1}\big[(b_1-b_2)(\F (t)+\Da \D{a}) -\I\big]^{-1}\W\nonumber \\
&~~~~~~~~-(\F (t)+ \Da \D{a})^{-1}\Z \F (t)^{-1}\nonumber \\
&~~~~~~~~-\W \G (t)^{-1}(\F (t)+ \Da \D{a})^{-1}-\W\Z\W \Big) \nonumber \\
&~~~~+b_2(b_1-b_2)\Big(- \G (t)^{-1}(\F (t)+ \Da \D{a})^{-1}\Z \nonumber \\
&~~~~~~~~- \big[(b_1-b_2)(\F (t)+\Da \D{a}) -\I\big]^{-1}\W \G (t)^{-1}\nonumber \\
&~~~~~~~~-\Z \F (t)^{-1}\big[(b_1-b_2)(\F (t)+\Da \D{a}) -\I\big]^{-1}\nonumber \\
&~~~~~~~~-\Z\W\Z \Big) \Big] \Bigg). \nonumber\
\end{align}

By the same argument used in the proof of Lemma \ref{O23.lem}, we observe that all components within the square brackets are element-wise negative and $(\Q{S_2}-\Q{S_1})$ is element-wise non-negative. Thus, $\frac{d}{dt} \hat{f}_a(\F (t)) \leq 0$.

%
%
%
Now, consider the following equality:
\begin{align*}
&\hat{f}_a(F(1)) = \\
&~~~~\hat{f}_a(F(0))+\int_0^1  \textstyle \frac{d}{dt}\hat{f}_a(\F (t))dt.
\end{align*}
Note that $\hat{f}_a(F(1)) = f_{a}(S_2)$ and $\hat{f}_a(F(0)) = f_{a}(S_1)$, and, as shown above,  $\frac{d}{dt}\hat{f_a}(\F (t))dt$ is non-positive. 
Thus,  $ f_{a}(S_1) \geq f_{a}(S_2)$. 

 From this we can conclude that $f_{a}$ is non-increasing for all $a \in V$, and, therefore, by Theorem~\ref{submod}, $f$ is submodular.

\balance

\bibliographystyle{IEEEtran}
\bibliography{thesis}


%
%



%

\end{document}